\documentclass[11pt]{amsart}
\usepackage[american]{babel}
\usepackage{graphicx}
\usepackage{amsmath}
\usepackage{amssymb}

\newtheorem{definition}{Definition}

\newtheorem{algo}{Algorithm}
\newtheorem{examplew}{Example}

\title{A GENERATING ALGORITHM FOR RIBBON TABLEAUX AND SPIN
  POLYNOMIALS} 
\author[Francois Descouens]{Francois Descouens} 
\address{Institut Gaspard
  Monge, Universit\'e de Marne-la-Vall\'ee\\ 77454 Marne-la-Vall\'ee
  Cedex 2, France}
\email{francois.descouens@univ-mlv.fr}

\begin{document}
\maketitle

\begin{abstract}
  We describe a general algorithm for generating various
  families of ribbon tableaux and computing their spin polynomials.
  This algorithm is derived from a new matricial coding. An advantage
  of this new notation lies in the fact that it permits one to
  generate ribbon tableaux with skew shapes.
\end{abstract}

\section{Introduction}
Ribbon tableaux are planar structures generalizing Young tableaux (see
\cite{Knuth, Lothaire} for the classical case). These are tilings of
Ferrers's diagram by {\it ribbons} (diagrams with special shape)
labelled with integers verifying some vertical and horizontal
monotonicity conditions.

Standard ribbon tableaux (all labels different) have first been
introduced by Stanton and White in 1985, in order to explain some
combinatorial properties of colored permutations \cite{SW}.
Semi-standard ribbon tableaux (repeated labels are allowed) go back to
the work of Lascoux, Leclerc and Thibon \cite{LLT}. These authors were
motivated by the introduction of $q$-analogues of certain
combinatorial identities, and in particular of the famous
Littlewood-Richardson rule describing products of Schur functions.
They obtained $q$-analogues of decomposition coefficients for any
product of Schur functions and many questions about these
$q$-coefficients are still open.

Studying ribbon tableaux is quite a difficult subject, which mainly
uses huge numerical experimentations. This is why we are interested in
finding efficient algorithms for generating and computing statistics
on them.

The matricial coding of ribbon tableaux used for numerical
experimentations in \cite{LLT,LT} is deficient because it does not
give some elementary properties (shape and position of the head of the
ribbons for example) without additional computations and cannot be
generalized to skew shapes. The algorithm used at the time was not
published and appeared only as a programming example with Maple/ACE in
\cite{KLLTUV} and a distributed version is described in \cite{Veign}.

The aim of this paper is to present a more general algorithm for
generating ribbon tableaux and computing spin polynomials, using a
different and more transparent coding for ribbon tableaux. These
algorithms are implemented in the combinatorial library {\tt
  MuPAD-Combinat}\cite{mu} which can be downloaded at
\textit{http://mupad-combinat.sourceforge.net/}.
\section{Basic definitions on ribbon tableaux}
We will mainly follow \cite{Andrews,Knuth,Lothaire} for classical
notions on partitions and tableaux and \cite{Lothaire} for notations
related to partitions. Let $\lambda$ and $\mu$ be two partitions such
that the diagram of $\lambda$ contains the diagram of $\mu$, The skew
partition of shape $\lambda / \mu$ can be defined as the set-theoretic
difference $\lambda-\mu$ ($\lambda$ is called the outer partition and
$\mu$ the inner).
\begin{figure}[h]
\begin{center}
\includegraphics*[width=4.5cm]{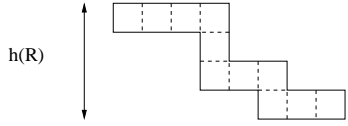}
\end{center}
\caption{An 11-ribbon of height $h(R)=4$}
\end{figure}
\begin{definition} Let $k$ be a nonnegative integer. A $k$-{\bf ribbon} 
  R is a connected skew diagram with $k$ cells which does not contains
  a 2$\times $2 square. The first (north-west) cell is called the {\bf
  head} and the last one (south-east) the {\bf tail}. The {\bf spin} is
  defined as $sp(R)=\frac{h(R)-1}{2}$.
\end{definition}
\begin{figure}[h]
\begin{center}
\includegraphics*[width=3cm]{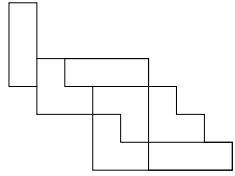}
\end{center}
\caption{A tiling with 3-ribbons of the skew partition $(8,7,6,5,1,1)
/ (3,3,1)$ of spin 3}
\end{figure}
We will denote by $Pav_k(\lambda / \mu)$ the set of $k$-ribbon tilings
of the skew shape $\lambda / \mu $. The {\bf spin} $sp(P)$ of a tiling
$P$ is the sum of the spins of its ribbons, and the {\bf cospin}
is the associated co-statistic into $Pav_k(\lambda / \mu)$, i.e:
\begin{displaymath}cosp(P)=\max \lbrace sp(U),~U\in Pav_k(\lambda / \mu) 
\rbrace ~-~sp(P).\end{displaymath} Ribbon tableaux are labelled 
tilings verifying monotonicity conditions similar to these of 
Young tableaux. The spin of a $k$-ribbon tableau
is the spin of the underlying tiling.
\begin{figure}[ht]
\begin{center}
\includegraphics*[width=3cm]{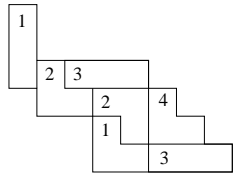}
\end{center}
\caption{A 3-ribbon tableau of shape $(8,7,6,5,1,1) / (3,3,1)$, 
weight $(2,2,2,1)$, and spin 3}
\end{figure}
\begin{definition}\index{ribbon tableaux!definition}
  A $k$-ribbon tableau of skew shape $\lambda / \mu$ is a tiling of
  the skew shape $\lambda /\mu$ by labelled $k$-ribbons such that the
  head of a ribbon labelled $i$ must not be on the right of a ribbon
  labelled $j>i$ and its tail must not be on the top of a ribbon
  labelled $j\ge i$. The weight of a $k$-ribbon tableau is the vector
  $\nu$ such that $\nu _i$ is the number of $k$-ribbons labelled $i$.
\end{definition}
We denote by $\text{Tab}_k(\lambda / \mu, \nu)$ the set of all
semi-standard $k$-ribbon tableaux of shape $\lambda / \mu$ and weight
$\nu$.
\begin{definition}\label{defpolyspin} The \textbf{spin} and 
\textbf{cospin polynomials} associated to the set
  $\text{Tab}_k(\lambda / \mu,\nu)$ are:
  \begin{displaymath}
  G_{\lambda / \mu, \nu}^{(k)}(q)=\sum_{T~\in
    \text{Tab}_k(\lambda / \mu, \nu)} q^{sp(T)}~~~\text{and}~~~
 \tilde{G}_{\lambda / \mu, \nu}^{(k)}(q)=\sum_{T~\in 
\text{Tab}_k(\lambda / \mu, \nu)} q^{cosp(T)}.
\end{displaymath} 
\end{definition}
 If we write $sp^{*}=\max \lbrace sp(T),~T\in Tab_k(\lambda / \mu)
\rbrace$, the following property holds:
\begin{displaymath}
 G_{\lambda / \mu, \nu}^{(k)}(q) = q^{sp^{*}}\tilde{G}_{\lambda / \mu,
 \nu}^{(k)}(\frac{1}{q}).
\end{displaymath}
\begin{examplew} In $Tab_3((8,7,6,5,1),(3,3,2,1))$, these two polynomials are:
\begin{displaymath}
G_{(8,7,6,5,1),(3,3,2,1)}^{(3)}(q)=3q^2+17q^3+33q^4+31q^5+18q^6+5q^7,
\end{displaymath}\begin{displaymath}
\tilde G_{(8,7,6,5,1),(3,3,2,1)}^{(3)}(q)= 3q^5+ 17q^4+33q^3+31q^2+18q
+5.
\end{displaymath}
\end{examplew}
\section{A generating algorithm}
In this section, we describe a new algorithm for generating all the
ribbon tableau in $\text{Tab}_k(\lambda / \mu, \delta)$ and computing
their spin polynomials. The main basic idea is to apply recursively
the algorithm of removing $k$-ribbon strips from a partition.
\subsection{A new coding for ribbon tableaux}
In this subsection, we extend the coding of \cite{Veign} which is not
well adapted to ribbon tableaux with skew shape $\lambda /
\mu$ and does not give immediate access to the shape of the tableau.
If the partition $\lambda$ is $(\lambda_1,\ldots,\lambda_n)$,
our coding is an $n \times \lambda_1$-array $(A_{i,j})$ defined as
follow:
\begin{itemize}
\item[1-] $A_{i,j}=-1$ if $(i,j) \in \mu$,
\item[2-] $A_{i,j}=p$ if there is the head of a ribbon labelled $p$ in the
  cell $(i,j)$,
\item[3-] $A_{i,j}=0$ otherwise.
\end{itemize}
By construction, we can immediately read the shape of the
corresponding ribbon tableau, and the length of the ribbons is
obtained by dividing the number of non negative cells by the number of
positive cells. This coding makes sense when $k=1$ because we obtain
the classical representation of a skew Young tableau (each ribbon is
reduced to its head).
\begin{examplew}
  The new coding of a 3-ribbon tableau of shape $(8,7,6,5,1,1)$ and weight $(3,3,2,1)$:\\
\begin{minipage}{5cm}
\begin{center}
\includegraphics*[width=3cm]{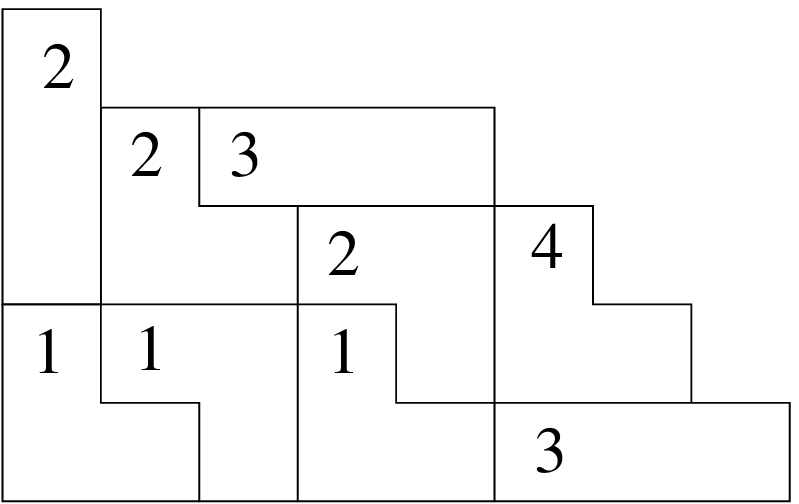}
\end{center}
\end{minipage}
\begin{minipage}{6cm}
\begin{center}
  $\longrightarrow~~~~~~~~~~~~~~~~~~~~~~~~é
\left ( \begin{array}{cccccccc}2& & & & & & & \\
                               0&2&3&0&0& & & \\
                               0&0&0&2&0&4& & \\
                               1&1&0&1&0&0&0& \\
                               0&0&0&0&0&3&0&0   \end{array}\right)$.
\end{center}
\end{minipage}
\end{examplew}
The structure used in \cite{Veign} was a matrix where the entry
$(i,j)$ is equal to $k$ if there is the tail of a ribbon labelled $i$
in the $j$-th column of the shape. We show that the main difficulty in
carrying out this construction is that we cannot read the shape of the
partition without reconstructing the entire ribbon tableau.
Furthermore, this coding become non injective when we try to
generalize it on ribbon tableaux of skew shape. With the previous
example we would obtain the following matrix:
\begin{displaymath}  
\left ( \begin{array}{cccccccc} 0&3&3&0&3&0&0&0\\ 3&0&3&0&3&0&0&0\\
                                  0&0&0&0&3&0&0&3\\
                                  0&0&0&0&0&0&3&0\end{array} \right).
\end{displaymath} 

\begin{examplew}
  The new coding of a 3-ribbon tableau of skew shape $(8,7,6,5,1,1) / (3,3,1)$ and weight $(3,2,2,1)$:\\
\vspace{0.5cm}
\begin{minipage}{5cm}
\begin{center}
\includegraphics*[width=3cm]{tabexemple.eps}
\end{center}
\end{minipage}
\begin{minipage}{6cm}
\begin{center}
  $\longrightarrow~~~~~~~~\left ( \begin{array}{cccccccc}1&  &  & & & & & \\
      0&  &  & & & & & \\
      0& 2& 3&0&0& & & \\
      -1& 0& 0&2&0&4& & \\
      -1&-1&-1&1&0&0&0& \\
      -1&-1&-1&0&0&3&0&0 \end{array}\right)$.
\end{center}
\end{minipage}
\end{examplew}
Decoding a $k$-ribbon tableau from an array $(A_{i,j})$
is as follow: for each label $p$ considered in decreasing order, we
fix vertically by the head a ribbon labelled $p$ in the cell $(i,j)$
of the shape if $A_{i,j}=p$, and we drive them to the frontier of the
partition.
\begin{algo}[From coding to ribbon tableaux]
$ $ \\ \vspace{-0.5cm}
\begin{itemize}
\item {\bf Input}: an array $M=(m_{ij})$.
\item {\tt Initialize}:
  \begin{itemize}
  \item $\lambda / \mu ~\longleftarrow$ {\tt the shape of the tableau coded by $M$},
   \item $\nu~\longleftarrow$ {\tt the
     weight of the tableau}. 
  \end{itemize}
 \item {\tt For each label $i$ in $\nu$:} 
     \begin{itemize}
     \item {\tt fix ribbons labelled $i$ by the head in each
         position corresponding to the cells labelled $i$ in $M$,}
     \item{\tt drive these ribbons on the frontier of $\lambda$,}
     \item$\lambda \longleftarrow \lambda$ {\tt without previous ribbons
     labelled $i$}.
      \end{itemize}
 \item {\bf Output}: The $k$-ribbon tableau corresponding to the array $M$.
\end{itemize}
\end{algo}
\begin{figure}[h] 
\begin{center}
  \includegraphics*[width=6.7cm]{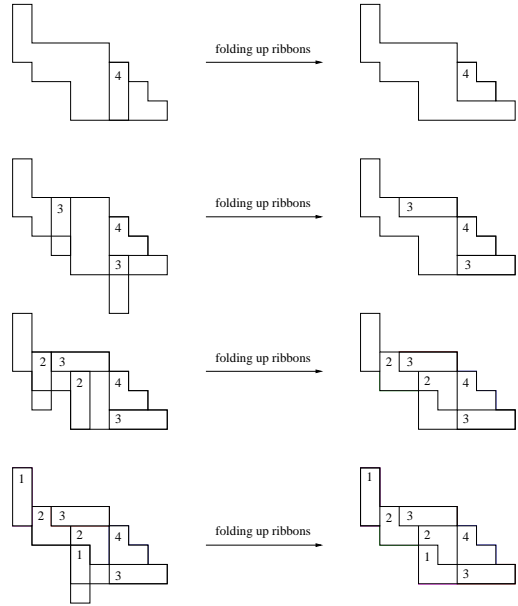}
\end{center}
\caption{Running of the algorithm on the previous 3-ribbon tableau of shape  
$(8,7,6,5,1,1) / (3,3,1)$ and weight $(2,2,2,1)$}
\end{figure}

\subsection{Adding and removing $k$-ribbon strip}
In order to explain our generating algorithm, we begin with a general
algorithm for adding or removing a $k$-ribbon strip from a partition
(these two operators on partitions come from the representation
theory of the quantum algebra $U_q(\widehat{sl_n})$, see
\cite{leclerc} for more details).
\begin{definition}
A skew tiling by $k$-ribbons $\Theta$ where the tail of each ribbon
is not on top of an other ribbon is called a $k$-{\bf ribbon
strip}. The {\bf weight} of $\Theta$ is the number of ribbons in the
tiling. Let $\Theta_\uparrow$ (resp. $\Theta_\downarrow$) be the the
horizontal strip made of the top cells (resp. the bottom cells) of the
columns of $\Theta$.
\end{definition}
\begin{figure}[h]
\begin{center}
\includegraphics[width=3cm, height=0.9cm]{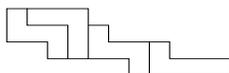}
\end{center}
\caption{A horizontal 5-ribbon strip of weight 4 and spin $3/2$}
\end{figure}
In a ribbon strip, the head of each ribbon lies in $\Theta_\uparrow$
and the tails in $\Theta_\downarrow$. As a $k$-ribbon strip has a
unique tiling by $k$-ribbons, it is completely defined by the
positions of all the ribbons's heads and the outer partition or all
the tails and the inner partition. For adding a $k$-ribbon strip
$\Theta$ to a partition, we represent this strip by the integer vector
$\mathcal{P}=(p_1,\ldots,p_n)$ with $p_i = k$ if $\Theta$ has a ribbon
with tail in the $i$-th column of $\Theta_{\downarrow}$. Similarly,
for removing a $k$-ribbon strip $\Theta$ from a partition, the vector
$\mathcal{P}$ is now defined by $p_i=-k$ if $\Theta$ has a ribbon with
head in the $i$-th column of $\Theta_{\uparrow}$.
\begin{algo}[Algorithm for adding a $k$-ribbon strip to a partition]
$ $ \\ \vspace{-0.5cm}
\begin{itemize}
\item {\bf Input}: the partition
  $\lambda=(\lambda_1,\ldots,\lambda_n)$ and the positions vector
  $\mathcal{P}=(p_1,\ldots,p_r)$.
\item {\tt Initialize}: 
\begin{itemize}
\item $\delta \longleftarrow (max(n,r)-1,\ldots,0)$,
\item $\lambda \longleftarrow$ {\tt conjugate partition of $\lambda$}.  
\end{itemize}
\item $\lambda \longleftarrow \lambda + \delta + \mathcal{P} $.
\item {\tt Sorting $\lambda$:} 
  \begin{itemize}
  \item $\lambda \longleftarrow \sigma (\lambda)$,
  \item $I \longleftarrow$ {\tt the inversions number of} $\sigma$.  
  \end{itemize}
\item {\tt if $\lambda - \delta$ is a partition:} 
  \begin{itemize}
  \item {\tt then }$\lambda \longleftarrow$ {\tt conjugate} $(\lambda -
  \delta)$ {\tt and} $I \longleftarrow \frac{(k-1)(p_1+\ldots +p_r) -
    I}{2k}$, 
\item {\tt else} $\lambda \longleftarrow$ {\tt FAIL}.
  \end{itemize}
\item {\bf Output}: $(\lambda, I)$ {\tt if exists, FAIL otherwise.}
\end{itemize} 
\end{algo}

In the previous algorithm, $\sigma$ corresponds to a permutation which
permits to sort the vector $\lambda$. Let defined the number of
inversions of the permutation $\sigma$ by the cardinality of the set
$\lbrace (i,j)\ \text{such that}\ i<j \ \text{and} \ \sigma(i) >
\sigma(j) \rbrace$. 
\subsection{Generating algorithm and computation of the spin and cospin 
polynomials for ribbon tableaux of a given shape and weight}
We will generalize the generating algorithm of \cite{Veign} to the
case of $k$-ribbon tableaux of skew shapes. A basic remark consists
in the fact that the ribbons labelled $i$ form a $k$ ribbon strip
$\Theta _i$ of weight $\nu_i$. That's why we search recursively all
the $k$-ribbon strip of weight $\nu _i$ contained in the shape and
removable from the outer partition and we fill an array with the
positions of the head in ${\Theta_i}{\uparrow}$. As the algorithm 1
also returns the spin of the added/removed $k$-ribbon strip, by
keeping the spin in each step of the construction we obtain finally
the spin of each tableau. In other words this algorithm permits one to
compute the spin and cospin polynomials without additional
computation.
\begin{algo}[Generating algorithm for ribbon tableaux and spin polynomials]
$ $ \\ \vspace{-0.5cm}
\begin{itemize}
\item {\bf Input}: The shape $\lambda / \mu$, the weight
  $\nu=(\nu_1,\ldots,\nu_p)$ and $k$.
\item {\tt Initialize:} $P=0$ {\tt and} $L\longleftarrow \lbrace (T,\lambda,0) \rbrace$
  {\tt where} $T$ {\tt is an array filled with -1 in cells
  corresponding to} $\mu$ {\tt and 0 otherwise.}
\item {\tt For each weight} $j$ {\tt from} $p$ {\tt down to 1},\\ 
      {\tt For each} $(T,\lambda_T,sp_T)$ {\tt in $L$ and each
permutation} $\mathcal{C}$ {\tt of the vector} $(\underbrace{0,\ldots
,0}_{\lambda_1-\nu_j},\underbrace{k,\ldots ,k}_{\nu_j~ times})$ {\tt
corresponding to a valid removed $k$-ribbon strip,}
      \begin{itemize}
         \item {\tt fill in $T$ the cells corresponding to the
frontier of $\lambda_T$ and the non-zero coordinates of $\mathcal{C}$,}
\item[] $\lambda_T \longleftarrow \lambda_T$
{\tt without the $k$-ribbon strip corresponding to $\mathcal{C}$,}
         \item $sp_T \longleftarrow sp_T +$ {\tt the spin of the previous $k$-ribbon strip,}
         \item {\tt if} $\mu \subset \lambda_T$ {\tt then} $L \longleftarrow L~\cup~(T,\lambda_T,sp_T).$
\end{itemize}
\item {\tt For each tableau $T$ in $L$, $P\longleftarrow P + q^{sp_T}$.} 
\item {\bf Output}: 
\begin{itemize}
\item The list $L$ of all the $k$-ribbon tableaux of skew shape $\lambda
  / \mu$ and weight $\nu$,
\item P the spin polynomial of this set. 
\end{itemize}
\end{itemize} 
\end{algo}
By studying the progress of this algorithm with huge set of ribbon
tableaux we remark that at the bottom of the tree there are a lot of
nodes, but in there is only few different shape for the remaining
partition. For example the number of nodes for $\lambda=(6,6,6,6,6,6)$
and $\nu=(3,1,1,1,1,1,1,1,1,1)$ is the sequence 3, 12, 48, 198, 780,
2940, 10080, 31080, 81480 and finally 43680 but the number of
different remaining shape are 3, 9, 16, 27, 33, 38, 33, 27, 16, 1. Our
previous algorithm search, at step $i$, all the possibility to retire
a $k$ ribbon strip from the remaining partition, but in fact we need
to test all the possibility only on the few remaining partitions.
That's why a recursive implementation of the previous algorithm with a
remember option seems to be the most efficient way to generate ribbon
tableaux.
\begin{examplew}
  In the case of 3-ribbons with $\lambda=(9,9,9,9,9,9,9,9,9)$ and
  standard weight $\mu=(1^{27})$, we have the following spin
  polynomial:\\\\
  \begin{minipage}{15cm}
  $
  414315330+8286306600q+85027356570q^2+588666753870q^3+3062543589300q^4+
  12659483135520q^5 + 42941179272810 q^6 + 121912682783970 q^7 +
  293410572110760 q^8 + 603798294330270 q^9 + 1068859924958280 q^{10}
  + 1634693172838050 q^{11} + 2166452577489720 q^{12} +
  2492870571244950 q^{13} + 2492870571244950 q^{14} + 2166452577489720
  q^{15} + 1634693172838050 q^{16} + 1068859924958280 q^{17} +
  603798294330270 q^{18} + 293410572110760 q^{19} + 121912682783970
  q^{20} + 42941179272810 q^{21} + 12659483135520 q^{22} +
  3062543589300 q^{23} + 588666753870 q^{24} + 85027356570 q^{25} +
  8286306600 q^{26}   + 414315330 q^{27}$
  \end{minipage}
  \\\\
  which correspond to a computation over $16\ 882\ 686\ 792\ 972\ 000$ ribbon
  tableaux.
\end{examplew}
\section*{Acknowledgments}
The author wishes to express his gratitude to the members of the
algebraic combinatorics team of the University of Marne-la-Vall\'ee
for their helpful comments.


\end{document}